**Yakiv O. Kalinovsky,** Dr.Sc., Senior Researcher,
Institute for Information Recording National Academy of Science of Ukraine, Kyiv, Shpaka str. 2, 03113, Ukraine, E-mail: kalinovsky@i.ua

**Dmitry V. Lande,** Dr.Sc., Head of Department,
Institute for Information Recording National Academy of Science of Ukraine, Kyiv, Shpaka str. 2, 03113, Ukraine, E-mail: dwlande@gmail.com

**Yuliya E. Boyarinova,** PhD, Associate Professor,
National Technical University of Ukraine "KPI", Kyiv, Peremogy av. 37, 03056, Ukraine, E-mail: ub@ua.fm

**Alina S. Turenko**, postgraduate student,
Institute for Information Recording National Academy of Science of Ukraine, Kyiv, Shpaka str. 2, 03113, Ukraine, E-mail: asturenko@mail.ru

# Computing Characteristics of One Class of Non-commutative Hypercomplex Number Systems of 4-dimension

The class of non-commutative hypercomplex number systems (HNS) of 4-dimension constructed by using of non-commutative procedure of Grassman-Clifford doubling of 2-dimensional systems is investigated in the article. All HNS of this class are constructed, algorithms of performance of operations and methods of algebraic characteristics calculation in them, such as conjugation, normalization, a type of zero dividers are investigated. Formulas of exponential functions representation in these systems are displayed.

## Introduction

The system of complex numbers which was opened in XVI century by Italian mathematicians D. Cardano and R. Bombeli, is widely used in science and technology and to this day is one of the main mathematical apparatuses for many theories and appendices. It, in turn, stimulated researches in search of greater dimensional number systems which would allow to solve so effectively not only problems on the plane, but also on three-dimensional and higher dimensions. For a long time such systems could not build. As it became clear later, these attempts were unsuccessful because researchers tried to receive systems with the same properties of real $\boldsymbol{R}$ and complex $\boldsymbol{C}$ numbers as commutativity, associativity and lack of zero divisors. By the way, it should be noted, there are not the same properties in complex numbers as in real ones - they can not be linearly arranged.

And only V. R. Hamilton managed to construct such number system of fourth dimension which possesses properties of division and associativity, but is non-commutative. It is about the quaternion system $\boldsymbol{H}$, for the first time presented by V.R. Hamilton in 1844 in work [1].

The quaternion system $\boldsymbol{H}$ found broad application in many scientific directions: in mechanics of a solid body for the rotation description in space, at the solution of problems of navigation, orientation and movement management; in computer animation, research of deformation of elastic designs, filtration of color images, cryptography and others. Which are possible to examine in work [2] in more detail.

Due to the above, it is advisable to investigate others non-commutative hypercomplex number systems (HNS) of fourth dimension. In principle, there are many such HNS. As to study all systems is very difficult, it is expedient to investigate separate classes of such HNS. Because quaternions $\boldsymbol{H}$ are a product of non-commutative doubling by means of Grassman-Clifford procedure of system of complex numbers $\boldsymbol{C}$ by the same system of complex numbers, authors decided to be limited by the class GChS of the fourth dimension, received by a way of non-commutative doubling of main HNS of second dimension by means of Grassman-Clifford procedure: the complex - $\boldsymbol{C}$, double - $\boldsymbol{W}$ and dual - $\boldsymbol{D}$ systems.

This work is devoted to research of arithmetic and algebraic operations and procedures in a class of such HNS and creation of exponential functions.

### Procedures of HNS doubling

Recurrent procedures of HNS doubling allow to build ranks of higher dimensions HNS. There are two types of procedures of doubling: Keli-Dixon procedure (CD-procedure) and Grassman-Clifford (GC-procedure) procedure.

The CD-procedure allows to build normed HNS of $2^n$-dimensions, where $n \in \boldsymbol{N}$ - an order of doubling [3–8]. If $n \geq 4$ we receive only nonassociative HNS.

GC-procedure allows to receive HNS which have more opportunities such as for dimension and properties [3, 5, 9].

We will consider the GC-procedure of doubling in more detail because it will be widely used further in this work.

We introduce following designations. In more general case we will designate of HNS as: $H(e,n)$, where $e = \{e_1,...,e_n\}$ – basis, and $n$ – dimension of HNS.

In the case, when we speak about concrete type of HNS, it will be denoted by name of its type as, for example, the system of complex numbers $\boldsymbol{C}\ (e)$. Here dimension may be omitted, as it is known from the type of HNS. But the identifier of basis $e$ should be given because doubling can consider two HNS of one type, but their bases should be distinguished among themselves.

We will designate doubling of system $H_1(e,m)$ by system $H_2(f,2)$ by using non-commutative GC procedure, as follows:

$$D(H_1(e,m), H_2(f,2)) = H_3(g,2m),$$

where $D$ – operator of doubling, and $2m$ – dimension of HNS $H_3$, which is the result of doubling. $2m$ basis elements $g$ will be all kinds of products of elements of bases $e$ and $f$:

$$g = \{e_1 f_1, e_1 f_2, e_2 f_1, ..., e_m f_2\}.$$

Cayley table consists of products of basis elements whose values reflect the specific properties of HNS.

### Definition of researched class of HNS

The class of non-commutative HNS of the fourth dimension, investigated in this work, consists of non-commutative doublings of HNS of the second dimension using GC procedure.

The basis of such HNS consists of four elements

$$g = \{g_1, g_2, g_3, g_4\} = \{e_1 f_1, e_1 f_2, e_2 f_1, e_2 f_2\}.$$

Certainly, the two-character elements can be replaced by one-character elements with indexes. But it will be made later as at this stage it is appropriate to use the two-character elements.

Investigated class of HNS will be determined by the following conditions:

1. Elements of bases $e_1$ and $f_1$ - unit elements of their systems;
2. Elements of Cayley table of HNS g are products in form: $g_i g_k = e_j f_s e_t f_r$, whose values can be calculated by switching multipliers and using Cayley tables of HNS which are doubling; thus we will consider that $e_1$ and $f_1$ commutative with $e_2$ and $f_2$, that is $e_1 f_2 = f_2 e_1, e_2 f_1 = f_1 e_2$, and the last are non-commutative among themselves, that is $e_2 f_2 = -f_2 e_2$.

For example:

$$g_1 g_1 = e_1 f_1 e_1 f_1 = e_1 e_1 f_1 f_1 = e_1 f_1 = g_1;$$
$$g_2 g_1 = e_1 f_2 e_1 f_1 = e_1 e_1 f_2 f_1 = e_1 f_2 = g_2;$$
$$g_2 g_3 = e_1 f_2 e_2 f_1 = -e_1 e_2 f_2 f_1 = -e_2 f_2 = -g_4;$$
$$g_2 g_2 = e_1 f_2 e_1 f_2 = e_1 e_1 f_2 f_2 = e_1 f_2 f_2;$$
$$g_4 g_4 = e_2 f_2 e_2 f_2 = -e_2 e_2 f_2 f_2.$$

The last two examples of elements of the Cayley table can be finished only for concrete HNS which are doubling.

3. Elements of the Cayley table, which are under the main diagonal, but not in the first column, are opposite for elements symmetric with respect to the main diagonal, that is

$$g_3 g_2 = -g_2 g_3; g_4 g_2 = -g_2 g_4; g_4 g_3 = -g_3 g_4.$$

Given these conditions, the generalized Cayley table for HNS of investigated class will be as follows:

| | $g_1$ | $g_2$ | $g_3$ | $g_4$ |
|---|---|---|---|---|
| $g_1$ | $g_1$ | $g_2$ | $g_3$ | $g_4$ |
| $g_2$ | $g_2$ | $e_1 f_2 f_2$ | $-g_4$ | $-e_2 f_2 f_2$ |
| $g_3$ | $g_3$ | $g_4$ | $e_2 e_2 f_1$ | $e_2 e_2 f_2$ |
| $g_4$ | $g_4$ | $e_2 f_2 f_2$ | $-e_2 e_2 f_2$ | $-e_2 e_2 f_2 f_2$ |

(1)

As it is known [3,9,10] there are three isomorphisms classes of HNS of second dimension. We will choose from these classes one representatives: the system of complex numbers $\boldsymbol{C}$, the system of double numbers $\boldsymbol{W}$ and the system of dual numbers $\boldsymbol{D}$.

As shown in the work [9], the first two operands in the operator of doubling can be commutative, as received Cayley tables differ only by order of lines and columns, that's why they are isomorphic.

Taking it into account, studied class HNS consists of six representatives classes of isomorphism:

1. $D(\boldsymbol{C}, \boldsymbol{C}, 4) = \boldsymbol{H}$ – quaternion system;
2. $D(\boldsymbol{C}, \boldsymbol{W}, 4) = \boldsymbol{AH}$ – antiquaternion system;
3. $D(\boldsymbol{C}, \boldsymbol{D}, 4) = D(\boldsymbol{C}, \boldsymbol{D}, 4)$;
4. $D(\boldsymbol{W}, \boldsymbol{W}, 4)$;
5. $D(\boldsymbol{D}, \boldsymbol{D}, 4)$;
6. $D(\boldsymbol{W}, \boldsymbol{D}, 4) = D(\boldsymbol{D}, \boldsymbol{W}, 4)$;

**Cayley tables of HNS of fourth dimension**

Cayley tables of the above six classes of isomorphisms can be easily obtained by substituting in (1) the basic elements of the complex – $\boldsymbol{C}$, double – $\boldsymbol{W}$ and dual numbers – $\boldsymbol{D}$, respectively

Having executed this algorithm for each of six cases, we will receive table 1

Table 1

**Cayley tables of hypercomplex number systems of fourth dimension**

| № | Designation | Cayley table |
|---|---|---|
| 1. | $\boldsymbol{H} = D(\boldsymbol{C}, \boldsymbol{C}, 4)$ | $\begin{array}{c\|cccc} H & e_1 & e_2 & e_3 & e_4 \\ \hline e_1 & e_1 & e_2 & e_3 & e_4 \\ e_2 & e_2 & -e_1 & e_4 & -e_3 \\ e_3 & e_3 & -e_4 & -e_1 & e_2 \\ e_4 & e_4 & e_3 & -e_2 & -e_1 \end{array}$ |
| 2. | $\boldsymbol{AH} = D(\boldsymbol{C}, \boldsymbol{W}, 4)$ | $\begin{array}{c\|cccc} AH & e_1 & e_2 & e_3 & e_4 \\ \hline e_1 & e_1 & e_2 & e_3 & e_4 \\ e_2 & e_2 & -e_1 & e_4 & -e_3 \\ e_3 & e_3 & -e_4 & e_1 & -e_2 \\ e_4 & e_4 & e_3 & e_2 & e_1 \end{array}$ |
| 3. | $D(\boldsymbol{C}, \boldsymbol{D}, 4)$ | $\begin{array}{c\|cccc} D(C,D,4) & e_1 & e_2 & e_3 & e_4 \\ \hline e_1 & e_1 & e_2 & e_3 & e_4 \\ e_2 & e_2 & -e_1 & e_4 & -e_3 \\ e_3 & e_3 & -e_4 & 0 & 0 \\ e_4 & e_4 & e_3 & 0 & 0 \end{array}$ |
| 4. | $D(\boldsymbol{W}, \boldsymbol{W}, 4)$ | $\begin{array}{c\|cccc} D(W,W,4) & e_1 & e_2 & e_3 & e_4 \\ \hline e_1 & e_1 & e_2 & e_3 & e_4 \\ e_2 & e_2 & e_1 & e_4 & e_3 \\ e_3 & e_3 & -e_4 & -e_1 & e_2 \\ e_4 & e_4 & -e_3 & -e_2 & e_1 \end{array}$ |
| 5. | $D(\boldsymbol{D}, \boldsymbol{D}, 4)$ | $\begin{array}{c\|cccc} D(D,D,4) & e_1 & e_2 & e_3 & e_4 \\ \hline e_1 & e_1 & e_2 & e_3 & e_4 \\ e_2 & e_2 & 0 & e_4 & 0 \\ e_3 & e_3 & -e_4 & 0 & 0 \\ e_4 & e_4 & 0 & 0 & 0 \end{array}$ |

| | | | |
|---|---|---|---|
| 6. | $D(\boldsymbol{W},\boldsymbol{D},4)$ | | $\begin{array}{c\|cccc} D(W,D,4) & e_1 & e_2 & e_3 & e_4 \\ \hline e_1 & e_1 & e_2 & e_3 & e_4 \\ e_2 & e_2 & e_1 & e_4 & e_3 \\ e_3 & e_3 & -e_4 & 0 & 0 \\ e_4 & e_4 & -e_3 & 0 & 0 \end{array}$ |

We will designate numbers of each of these systems as.

$$w = a_1e_1 + a_2e_2 + a_3e_3 + a_4e_4, \tag{2}$$

where: $a_i \in R$.

### Addition and multiplication operations

In these systems are entered addition and multiplication operations as follows:

The number $w_3$: $w_3 = w_1 + w_2 = (a_1 + b_1)e_1 + (a_2 + b_2)e_2 + (a_3 + b_3)e_3 + (a_4 + b_4)e_4$ is called *the sum* of two numbers $w_1 = a_1e_1 + a_2e_2 + a_3e_3 + a_4e_4$ and $w_2 = b_1e_1 + b_2e_2 + b_3e_3 + b_4e_4$ for all systems from table 1.

As for the product, for each of systems it will have other appearance because depends on Cayley table in considered hypercomplex number system. The product of two elements of each of the systems we present in Table 2.

Table 2.

**Rules of multiplication**

| **№** | **HNS** | **Rule of multiplication** |
|---|---|---|
| 1. | $\boldsymbol{H}$ | $w_1w_2 = (a_1b_1 - a_2b_2 - a_3b_3 - a_4b_4)e_1 + (a_1b_2 + a_2b_1 + a_3b_4 - a_4b_3)e_2 + (a_1b_3 + a_3b_1 - a_2b_4 + a_4b_2)e_3 + (a_1b_4 + a_4b_1 + a_2b_3 - a_3b_2)e_4$ |
| 2. | $\boldsymbol{AH}$ | $w_1w_2 = (a_1b_1 - a_2b_2 + a_3b_3 + a_4b_4)e_1 + (a_1b_2 + a_2b_1 - a_3b_4 + a_4b_3)e_2 + (a_1b_3 + a_3b_1 - a_2b_4 + a_4b_2)e_3 + (a_1b_4 + a_4b_1 + a_2b_3 - a_3b_2)e_4$ |
| 3. | $D(\boldsymbol{C},\boldsymbol{D},4)$ | $w_1w_2 = (a_1b_1 - a_2b_2)e_1 + (a_1b_2 + a_2b_1)e_2 + (a_1b_3 + a_3b_1 - a_2b_4 + a_4b_2)e_3 + (a_1b_4 + a_4b_1 + a_2b_3 - a_3b_2)e_4$ |
| 4. | $D(\boldsymbol{W},\boldsymbol{W},4)$ | $w_1w_2 = (a_1b_1 + a_2b_2 + a_3b_3 - a_4b_4)e_1 + (a_1b_2 + a_2b_1 - a_3b_4 + a_4b_3)e_2 + (a_1b_3 + a_3b_1 + a_2b_4 - a_4b_2)e_3 + (a_1b_4 + a_4b_1 + a_2b_3 - a_3b_2)e_4$ |
| 5. | $D(\boldsymbol{D},\boldsymbol{D},4)$ | $w_1w_2 = a_1b_1e_1 + (a_1b_2 + a_2b_1)e_2 + (a_1b_3 + a_3b_1)e_3 + (a_1b_4 + a_4b_1 + a_2b_3 - a_3b_2)e_4$ |
| 6. | $D(\boldsymbol{W},\boldsymbol{D},4)$ | $w_1w_2 = (a_1b_1 + a_2b_2)e_1 + (a_1b_2 + a_2b_1)e_2 + (a_1b_3 + a_3b_1 + a_2b_4 - a_4b_2)e_3 + (a_1b_4 + a_4b_1 + a_2b_3 - a_3b_2)e_4$ |

According to rules of addition and multiplication it is possible to mark out their main properties:

1) addition operation is commutative:

$$w_1 + w_2 = w_2 + w_1;$$

2) addition operation is associative:

$$(w_1 + w_2) + w_3 = w_1 + (w_2 + w_3);$$

3) multiplication operation is noncommutative:

$$w_1 w_2 \neq w_2 w_1. \qquad (3)$$

We will prove, for example, for system $D(\boldsymbol{C}, \boldsymbol{D}, 4)$, really:

$$\begin{aligned} w_1 w_2 &= (a_1 e_1 + a_2 e_2 + a_3 e_3 + a_4 e_4)(b_1 e_1 + b_2 e_2 + b_3 e_3 + b_4 e_4) = \\ &= (a_1 b_1 - a_2 b_2) e_1 + (a_1 b_2 + a_2 b_1) e_2 + \\ &+ (a_1 b_3 - a_2 b_4 + a_3 b_1 + a_4 b_2) e_3 + (a_1 b_4 + a_2 b_3 - a_3 b_2 + a_4 b_1) e_4 \end{aligned},$$

but opposite order is such as:

$$\begin{aligned} w_2 w_1 &= (b_1 e_1 + b_2 e_2 + b_3 e_3 + b_4 e_4)(a_1 e_1 + a_2 e_2 + a_3 e_3 + a_4 e_4) = \\ &= (b_1 a_1 - b_2 a_2) e_1 + (b_1 a_2 + b_2 a_1) e_2 + \\ &+ (b_1 a_3 - b_2 a_4 + b_3 a_1 + b_4 a_2) e_3 + (b_1 a_4 + b_2 a_3 - b_3 a_2 + b_4 a_1) e_4 = \\ &= (a_1 b_1 - a_2 b_2) e_1 + (a_1 b_2 + a_2 b_1) e_2 + \\ &+ (a_1 b_3 + a_2 b_4 + a_3 b_1 - a_4 b_2) e_3 + (a_1 b_4 - a_2 b_3 + a_3 b_2 + a_4 b_1) e_4 \neq w_1 w_2 \end{aligned}$$

That is carried out (3).

For other systems the proof is similar.

4) multiplication operation is associative:

$$w_1 (w_2 w_3) = (w_1 w_2) w_3.$$

It can be proved directly, using table 2.

5) In the same way it is possible to prove the distributivity:

$$w_1 (w_2 + w_3) = w_1 w_2 + w_1 w_3;$$

6) for each system is determined action of multiplication by a scalar $k \in R$,

$$k w_1 = k a_1 e_1 + k a_2 e_2 + k a_3 e_3 + k a_4 e_4;$$

7) for $\forall k_1, k_2 \in R$ is performed $(k_1 w_1)(k_2 w_2) = k_1 k_2 (w_1 w_2)$.

**Norm definition**

In the work [12] the norm of hypercomplex number generally is determined by a formula

$$N(w) = \sum_{i=1}^{n} \gamma_{ij}^{k} a_i, \qquad (4)$$

where $\gamma_{ij}^{k}$ - structural constants of hypercomplex number system of, which are defined from table 1. On this basis the norm of matrix is constructed [12].

For system $D(W, W, 4)$ the matrix of norm will have an appearance

$$N(w)=\begin{vmatrix} a_1 & a_2 & a_3 & -a_4 \\ a_2 & a_1 & a_4 & -a_3 \\ a_3 & -a_4 & a_1 & a_2 \\ a_4 & -a_3 & a_2 & a_1 \end{vmatrix}. \tag{5}$$

Having calculated the determinant of matrix (8) we will receive norm of hypercomplex number $w$:

$$N(w)=\left(a_1{}^2-a_2{}^2-a_3{}^2+a_4{}^2\right)^2. \tag{6}$$

By analogy to the theory of quaternions we will call a root of the norm a pseudonorm of hypercomplex number (6), which will be denoted as $N(w)$:

$$N(w)=a_1{}^2-a_2{}^2-a_3{}^2+a_4{}^2. \tag{7}$$

For other hypercomplex number systems the norm is determined by the same algorithm. It should be noted that for each of these systems the matrix of norm will have other appearance, and according to it is differ pseudonorm representations that it is possible to see from table 3.

Table 3

**Pseudonorm**

| № | HNS | Matrix of norm | Pseudonorm |
|---|---|---|---|
| 1. | $\boldsymbol{H}$ | $N(w)=\begin{vmatrix} a_1 & -a_2 & -a_3 & -a_4 \\ a_2 & a_1 & a_4 & -a_3 \\ a_3 & a_4 & a_1 & -a_2 \\ a_4 & -a_3 & a_2 & a_1 \end{vmatrix}$ | $N(w)=a_1^{\;2}+a_2^{\;2}+a_3^{\;2}+a_4^{\;2}$ |
| 2. | $\boldsymbol{AH}$ | $N(w)=\begin{vmatrix} a_1 & -a_2 & a_3 & a_4 \\ a_2 & a_1 & a_4 & -a_3 \\ a_3 & a_4 & a_1 & -a_2 \\ a_4 & -a_3 & a_2 & a_1 \end{vmatrix}$ | $N(w)=a_1^{\;2}+a_2^{\;2}-a_3^{\;2}-a_4^{\;2}$ |
| 3. | $D(\boldsymbol{C},\boldsymbol{D},4)$ | $N(w)=\begin{vmatrix} a_1 & -a_2 & 0 & 0 \\ a_2 & a_1 & 0 & 0 \\ a_3 & a_4 & a_1 & -a_2 \\ a_4 & -a_3 & a_2 & a_1 \end{vmatrix}$ | $N(w)=a_1^{\;2}+a_2^{\;2}$ |
| 4. | $D(\boldsymbol{W},\boldsymbol{W},4)$ | $N(w)=\begin{vmatrix} a_1 & a_2 & a_3 & -a_4 \\ a_2 & a_1 & a_4 & -a_3 \\ a_3 & -a_4 & a_1 & a_2 \\ a_4 & -a_3 & a_2 & a_1 \end{vmatrix}$ | $N(w)=a_1^{\;2}-a_2^{\;2}-a_3^{\;2}+a_4^{\;2}$ |
| 5. | $D(\boldsymbol{D},\boldsymbol{D},4)$ | $N(w)=\begin{vmatrix} a_1 & 0 & 0 & 0 \\ a_2 & a_1 & 0 & 0 \\ a_3 & 0 & a_1 & 0 \\ a_4 & -a_3 & a_2 & a_1 \end{vmatrix}$ | $N(w)=a_1^{\;2}$ |
| 6. | $D(\boldsymbol{W},\boldsymbol{D},4)$ | $N(w)=\begin{vmatrix} a_1 & a_2 & 0 & 0 \\ a_2 & a_1 & 0 & 0 \\ a_3 & -a_4 & a_1 & a_2 \\ a_4 & -a_3 & a_2 & a_1 \end{vmatrix}$ | $N(w)=a_1^{\;2}-a_2^{\;2}$ |

Apparently from table 3, in some systems the pseudonorm can be negative. It is possible to show that the pseudonorm entered by such method is multiplicative for each of considered systems, i.e. equality is carried out:

$$N(w_1 w_2) = N(w_1) N(w_2). \quad (8)$$

**Definition and characteristics of conjugate numbers**

Let we have any number $w = a_1 e_1 + a_2 e_2 + a_3 e_3 + a_4 e_4$, conjugate to it we designate $\overline{w} = b_1 e_1 + b_2 e_2 + b_3 e_3 + b_4 e_4$, where $b_1, b_2, b_3, b_4$ – unknown coefficients.

As is offered in [12], definition of conjugate is entered on the basis of equality

$$w\overline{w} = N(w), \quad (9)$$

Substituting introduced notation and using Table 2, and equating coefficients of the same basic elements we obtain a linear algebraic system with respect to variables $b_1, b_2, b_3, b_4$.

For hypercomplex number system $D(\boldsymbol{D}, \boldsymbol{D}, 4)$, for example, this linear algebraic system is

$$\begin{cases} a_1 b_1 = a_1^{\ 2} \\ a_1 b_2 + a_2 b_1 = 0 \\ a_1 b_3 + a_3 b_1 = 0 \\ a_1 b_4 + a_4 b_1 + a_2 b_3 - a_3 b_2 = 0 \end{cases}, \quad (10)$$

which solutions have the form:

$$b_1 = a_1, b_2 = -a_2, b_3 = -a_3, b_4 = -a_4. \quad (11)$$

Therefore, if the original number $w = a_1 e_1 + a_2 e_2 + a_3 e_3 + a_4 e_4$, that the conjugate number to it has view:

$$\overline{w} = a_1 e_1 - a_2 e_2 - a_3 e_3 - a_4 e_4. \quad (12)$$

In spite of the fact that for each of the systems given in table 1, the linear algebraic system (10) will have other appearance, representation of the interfaced number for each case will have an appearance (12).

We will define some properties of conjugate numbers.

1) the sum and the product of conjugate numbers are real numbers;

2) the conjugate of the sum is the sum of conjugated $\overline{w_1 + w_2} = \overline{w_1} + \overline{w_2}$;

3) the conjugate of the product is the product of conjugated $\overline{w_1 w_2} = \overline{w_1}\,\overline{w_2}$, which can be verified directly.

**Zero divisors and their properties**

Not equal to zero hypercomplex number $w_1 \neq 0$ is called as zero divider if there is such hypercomplex number $w_2 \neq 0$, that their product is equal to zero $w_1 w_2 = 0$, and it means the same ratio between their pseudonorm:

$$N(w_1 w_2) = 0. \quad (13)$$

On the basis of (8) pseudonorm of a divider of zero is equal to zero

$$N(w_1) = 0. \quad (14)$$

From (14) follows the signs of zero divisor in any hypercomplex number system, considered systems (except system of quaternions for which according to Frobenius's theorem there are no zero divisors), which we will give in table 4

Table 4

**Signs of zero divisor**

| № | HNS | Sign of zero divisor |
|---|---|---|
| 1. | $\boldsymbol{H}$ | – |
| 2. | $\boldsymbol{AH}$ | $a_1^2 + a_2^2 = a_3^2 + a_4^2$ |
| 3. | $D(\boldsymbol{C}, \boldsymbol{D}, 4)$ | $a_1 = a_2 = 0$ |
| 4. | $D(\boldsymbol{W}, \boldsymbol{W}, 4)$ | $a_1^2 + a_4^2 = a_3^2 + a_2^2$ |
| 5. | $D(\boldsymbol{D}, \boldsymbol{D}, 4)$ | $a_1 = 0$ |
| 6. | $D(\boldsymbol{W}, \boldsymbol{D}, 4)$ | $a_1 = \pm a_2$ |

**Division operation**

The share from the left division of hypercomplex number $w_1$ into hypercomplex number $w_2$ is the solution of the equation

$$w_2 x = w_1 \qquad (15)$$

To solve the equation (15) it is necessary to multiply its both parts at the left at first on $\overline{w_2}$, and then on $\frac{1}{|w_2|^2}$, where $|w_2|^2 \neq 0$. We will receive

$$x_l = \frac{1}{|w_2|^2}\overline{w_2} w_1 . \qquad (16)$$

Direct substitution of (16) in the equation (15) we find out that this expression is the solution of this equation.

We introduce the right division on the basis the equation

$$x w_2 = w_1 , \qquad (17)$$

from where

$$x_r = \frac{1}{|w_2|^2} w_1 \overline{w_2} . \qquad (18)$$

Since, we review non-commutatuve systems that the product of numbers depends on the order of the factors, then $x_l \neq x_r$. Thus, the solution of the equation (15) is called the left share, and the equation (17) - the right share [9].

It should be noted that the operation of division, unlike the fields of real and complex numbers, not possible, not only by zero but also by zero divisors in this systems.

**The principle of constructing exponential representation of hypercomplex variable by using of the associated system of differential equations**

This method consists in the following [12].

Representation of exponent in the system $H(e,n)$ of numbers $M \in H(e,n)$ is a particular solution of hypercomplex linear differential

$$\dot{X} = MX \,, \tag{19}$$

with the initial condition

$$Exp(0) = e_1 \,. \tag{20}$$

For creation of the solution of hypercomplex linear differential equation (19) we will present it in a vector-matrix form. Thus

$$\dot{\overline{X}} = (\dot{x}_1, ..., \dot{x}_n)^r \,, \tag{21}$$

and the vector column $\overline{MX}$ received from hypercomplex numbers $MX$, it is possible to present in the matrix form of some matrix $M$ of the dimension $n \times n$ which elements are linear combinations of component hypercomplex numbers $M$, on a vector column $\overline{X}$:

$$\overline{MX} = M\overline{X} \,. \tag{22}$$

Then the hypercomplex equations (19) can be presented in the form of the associated system of $n$ linear differential equations of the first order:

$$\dot{\overline{X}} = M\overline{X} \,. \tag{23}$$

Further it is necessary to find characteristic values $\lambda_1, ..., \lambda_{\ddot{\imath}}$ of a matrix $M$, that is to solve the characteristic equation

$$M - \lambda E = 0 \,. \tag{24}$$

Thus, characteristic values $\lambda_1, ..., \lambda_{\ddot{\imath}}$ will depend on hypercomplex number $M$.

Further it is necessary to construct the common decision which depends on $n^2$ arbitrary constants, from them $n^2 - \ddot{\imath}$ are linear and dependent on $n$ of arbitrary variables. For obtaining these linear dependences it is necessary to solve the system of linear equations then it is possible to receive the common decision (21) which depends on $n$ arbitrary constants - $\overline{X}(t, C_1, ..., C_n)$. Value of arbitrary constants are establish by using an entry condition (20). Components of decision vector column $\overline{X}$ also will be components of exponent of hypercomplex number $M$:

$$Exp(M) = \sum_{i=1}^{n} \overline{x}_i e_i \,. \tag{25}$$

At creation of representation of exponent for each of considered hypercomplex number systems it was revealed that in some of them exponent representation includes three cases, and in some - one. We will give examples.

**Construction of exponent representation of variable from system $D(W, W, 4)$**

According to the algorithm given in the previous point, we will consider four-dimensional number system $D(W, W, 4)$. We will build exponent representations of number $M \in D(W, W, 4)$ in it.

The equation (19) can be rewritten in a form

$$\frac{dX}{dt} = MX \,. \tag{26}$$

We calculate the right side of equation (26). According to the multiplication table for these numbers we have:

$$M = \left(m_1 e_1 + m_2 e_2 + m_3 e_3 + m_4 e_4\right),\ X = \left(x_1 e_1 + x_2 e_2 + x_3 e_3 + x_4 e_4\right)$$

$$\begin{aligned} MX &= \left(m_1 e_1 + m_2 e_2 + m_3 e_3 + m_4 e_4\right)\left(x_1 e_1 + x_2 e_2 + x_3 e_3 + x_4 e_4\right) = \\ &= \left(m_1 x_1 m_2 x_2 + m_3 x_3 - m_4 x_4\right)e_1 + \left(m_1 x_2 + m_2 x_1 - m_3 x_4 + m_4 x_3\right)e_2 + , \\ &+ \left(m_1 x_3 + m_2 x_4 + m_3 x_1 - m_4 x_2\right)e_3 + \left(m_1 x_4 + m_2 x_3 - m_3 x_2 + m_4 x_1\right)e_4 \end{aligned} \quad (27)$$

According to (27), equation (26) will have the form:

$$\begin{aligned} \frac{dX}{dt} &= \left(m_1 x_1 m_2 x_2 + m_3 x_3 - m_4 x_4\right)e_1 + \left(m_1 x_2 + m_2 x_1 - m_3 x_4 + m_4 x_3\right)e_2 + \\ &+ \left(m_1 x_3 + m_2 x_4 + m_3 x_1 - m_4 x_2\right)e_3 + \left(m_1 x_4 + m_2 x_3 - m_3 x_2 + m_4 x_1\right)e_4 \end{aligned} \quad (28)$$

From here it is possible to write down the associated system:

$$\begin{cases} \dfrac{dx_1}{dt} = m_1 x_1 + m_2 x_2 + m_3 x_3 - m_4 x_4 \\ \dfrac{dx_2}{dt} = m_2 x_1 + m_1 x_2 + m_4 x_3 - m_3 x_4 \\ \dfrac{dx_3}{dt} = m_3 x_1 - m_4 x_2 + m_1 x_3 + m_2 x_4 \\ \dfrac{dx_4}{dt} = m_4 x_1 - m_3 x_2 + m_2 x_3 + m_1 x_4 \end{cases} \quad (29)$$

For solution of the equation (29) it is necessary to find roots of the characteristic equation (24),

where $M = \begin{pmatrix} m_1 & m_2 & m_3 & m_4 \\ m_2 & m_1 & m_4 & -m_3 \\ m_3 & -m_4 & m_1 & m_2 \\ m_4 & -m_3 & m_2 & m_1 \end{pmatrix}$, $E$ – single matrix.

Taking into account the aforesaid, the characteristic equation will have the form:

$$\begin{vmatrix} m_1 - \lambda & m_2 & m_3 & m_4 \\ m_2 & m_1 - \lambda & m_4 & -m_3 \\ m_3 & -m_4 & m_1 - \lambda & m_2 \\ m_4 & -m_3 & m_2 & m_1 - \lambda \end{vmatrix} = 0 . \quad (30)$$

Having calculated the equation (30), we will receive the following equation:

$$\left((m_1 - \lambda)^2 - m_2^{\ 2} - m_3^{\ 2} + m_4^{\ 2}\right)^2 = 0 . \quad (31)$$

With (31) it is easy to find characteristic values

$$\lambda_{1,2,3,4} = m_1 \pm \sqrt{m_2^{\ 2} + m_3^{\ 2} - m_4^{\ 2}} , \quad (32)$$

as we see, they can be different depending on what sign has a radical expression. We will consider special cases:

1) $m_2^{\ 2} + m_3^{\ 2} - m_4^{\ 2} < 0$

We will introduce designations:

$$\overline{m} = \sqrt{m_2^{\ 2} + m_3^{\ 2} - m_4^{\ 2}} \ . \tag{33}$$

In this case, we will have multiple complex-conjugate roots of the characteristic equation, ie,:

$$\lambda_{1,2,3,4} = m_1 \pm i\sqrt{m_4^{\ 2} - m_2^{\ 2} - m_3^{\ 2}} \ .$$

For convenience we will introduce one more designation:

$$\overline{\overline{m}} = \sqrt{m_4^{\ 2} - m_2^{\ 2} - m_3^{\ 2}} \ . \tag{34}$$

Then roots of the characteristic equation (30) will have the form:

$$\lambda_{1,2,3,4} = m_1 \pm i\overline{m} \, . \tag{35}$$

According to (35) solution of system (29) will sought in the following form:

$$X_k(t) = e^{m_1 t}\left((A_k + B_k t)\cos\left(\overline{\overline{m}}t\right) + (C_k + D_k t)\cos\left(\overline{\overline{m}}t\right)\right), \ k = 1,\dots,4 \, . \tag{36}$$

We calculate the first derivative

$$\begin{aligned} \frac{dX_k(t)}{dt} = e^{m_1 t}\Big(&m_1 A_k + B_k + \overline{\overline{m}}C_k + \left(m_1 B_k + \overline{\overline{m}}D_k\right)t\cos\left(\overline{\overline{m}}t\right) + \\ &+ m_1 C_k + D_k - \overline{\overline{m}}A_k + \left(m_1 D_k + \overline{\overline{m}}B_k\right)t\sin\left(\overline{\overline{m}}t\right), \ k = 1,\dots,4 \end{aligned} \tag{37}$$

If we substitute equalities (36) and (37) in the associated system (29) and use a method of indeterminate coefficients we will receive the system of 16 equations with 16 unknown which can be reduced to a vector-matrix system by introducing the following designations::

$$Q = \begin{vmatrix} 0 & m_2 & m_3 & -m_4 \\ m_2 & 0 & m_4 & -m_3 \\ m_3 & -m_4 & 0 & m_2 \\ m_4 & -m_3 & m_2 & 0 \end{vmatrix}, \quad \overline{A} = \begin{vmatrix} A_1 \\ A_2 \\ A_3 \\ A_4 \end{vmatrix}, \quad \overline{B} = \begin{vmatrix} B_1 \\ B_2 \\ B_3 \\ B_4 \end{vmatrix}, \quad \overline{C} = \begin{vmatrix} C_1 \\ C_2 \\ C_3 \\ C_4 \end{vmatrix}, \quad \overline{D} = \begin{vmatrix} D_1 \\ D_2 \\ D_3 \\ D_4 \end{vmatrix}$$

The matrix $Q$ possesses such property:

$$QQ = -\overline{\overline{m}}^{\,2} E \, , \tag{38}$$

where $E$ – single matrix.

Having used such designations, the vector-matrix system will have the form:

$$\begin{cases} \overline{B} + \overline{\overline{m}}\overline{C} = Q\overline{A} \\ \overline{\overline{m}}\overline{D} = Q\overline{B} \\ -\overline{\overline{m}}\overline{A} + \overline{D} = Q\overline{C} \\ -\overline{\overline{m}}\overline{B} = Q\overline{D} \end{cases} \tag{39}$$

If the 4th equation is multiplied by a matrix *Q*, then $-\overline{\overline{m}}Q\overline{B} = QQ\overline{D}$, if to use a relation (38) and cut back on $\overline{\overline{m}}$, we will receive the 2nd equation. It means that the 2nd and 4th equation are linearly dependent, of the 2nd equation we will express $\overline{D}$ through $\overline{B}$: $\overline{D} = \frac{1}{\overline{\overline{m}}}Q\overline{B}$ and will substitute it in the 3rd, we will multiply the 1st equation at the left by *Q*. Taking into account (43) we will receive the following system:

$$\begin{cases} Q\overline{B} + \overline{\overline{m}}Q\overline{C} = QQ\overline{A} \\ -\overline{\overline{m}}\overline{A} + \frac{1}{\overline{\overline{m}}}\overline{B}Q = Q\overline{C} \end{cases}$$

The sum of the equations of this system will give us the equation $2Q\overline{B} = 0$ from which follows ,

$$\overline{B} = \overline{D} = 0 \tag{40}$$

that is

$$\overline{B} = \begin{vmatrix} 0 \\ 0 \\ 0 \\ 0 \end{vmatrix}, \quad \overline{D} = \begin{vmatrix} 0 \\ 0 \\ 0 \\ 0 \end{vmatrix}.$$

Then, if we substitute (40) in the 1st equation of system (39), we will receive

$$\overline{A} \text{ – arbitrary constant, } C = \frac{1}{\overline{\overline{m}}}Q\overline{A}. \tag{41}$$

Taking into account (40) and (41), (36) will have a vector-matrix appearance:

$$X(t) = e^{m_1 t}\left(\overline{A}\cos\left(\overline{\overline{m}}t\right) + \frac{1}{\overline{\overline{m}}}Q\overline{A}\sin\left(\overline{\overline{m}}t\right)\right). \tag{42}$$

This expression includes an arbitrary column vector consisting of 4 constants of integration. For their determination we use the initial condition, i.e. we will consider a reference point of *t=0*, *K=0.*

Substitution of these values in a row for exponent gives

$$Exp(0) = 1 \cdot e_1 + 0 \cdot e_2 + 0 \cdot e_3 + 0 \cdot e_4, \tag{43}$$

that is, must be:

$$\overline{X}(0)=\begin{vmatrix}1\\0\\0\\0\end{vmatrix}. \tag{44}$$

Substituting (44) into (36), we have:

$$A_1=1,\ A_2=0,\ A_3=0,\ A_4=0. \tag{45}$$

Solution of the system can be obtained if to substitute values $A_1$, $A_2$, $A_3$, $A_4$ into (36):

$$\begin{cases} x_1=e^{m_1t}\cos\left(\overline{\overline{m}}t\right) \\ x_2=e^{m_1t}\dfrac{m_2}{\overline{\overline{m}}}\sin\left(\overline{\overline{m}}t\right) \\ x_3=e^{m_1t}\dfrac{m_3}{\overline{\overline{m}}}\sin\left(\overline{\overline{m}}t\right) \\ x_4=e^{m_1t}\dfrac{m_4}{\overline{\overline{m}}}\sin\left(\overline{\overline{m}}t\right) \end{cases}. \tag{46}$$

By means of the common solutions it is possible to write down exponential function of variable from system $D(\boldsymbol{W},\boldsymbol{W},4)$, if to accept $m_it\approx m_i$

$$Exp(M)=e^{m_1}\left(\cos\overline{\overline{m}}\cdot e_1+\frac{\sin\overline{\overline{m}}}{\overline{\overline{m}}}\left(m_2e_2+m_3e_3+m_4e_4\right)\right)=$$

$$=e^{m_1}\left(\cos i\sqrt{\left|m_2^{\ 2}+m_3^{\ 2}-m_4^{\ 2}\right|}e_1-\frac{i\sin i\sqrt{\left|m_2^{\ 2}+m_3^{\ 2}-m_4^{\ 2}\right|}}{\sqrt{\left|m_2^{\ 2}+m_3^{\ 2}-m_4^{\ 2}\right|}}\left(m_2e_2+m_3e_3+m_4e_4\right)\right)$$

Using communications of trigonometrical functions with the hyperbolic we can write down

$$Exp(M)=e^{m_1}\left(ch\sqrt{\left|m_2^{\ 2}+m_3^{\ 2}-m_4^{\ 2}\right|}e_1+\frac{\left(m_2e_2+m_3e_3+m_4e_4\right)}{\sqrt{\left|m_2^{\ 2}+m_3^{\ 2}-m_4^{\ 2}\right|}}sh\sqrt{\left|m_2^{\ 2}+m_3^{\ 2}-m_4^{\ 2}\right|}\right). \tag{47}$$

2) $m_2^{\ 2}+m_3^{\ 2}-m_4^{\ 2}>0$

In this case we will have two multiple valid roots:

$$\lambda_{1,2,3,4}=m_1\pm\overline{m}. \tag{48}$$

The solution of system (29) is looked for in such form:

$$\begin{aligned} X_k(t)&=e^{\left(m_1+\overline{m}\right)t}\left(A_k+B_kt\right)+e^{\left(m_1-\overline{m}\right)t}\left(C_k+D_kt\right)= \\ &=e^{m_1t}\left(e^{\overline{m}t}\left(A_k+B_kt\right)+e^{-\overline{m}t}\left(C_k+D_kt\right)\right),\ k=1,\dots,4 \end{aligned} \tag{49}$$

We calculate the first derivative:

$$\frac{dX_k(t)}{dt} = e^{m_1 t}\left(e^{\overline{m}t}\left(\left(A_k + B_k t\right)\left(m_1 + \overline{m}\right) + B_k\right) + e^{-\overline{m}t}\left(\left(C_k + D_k t\right)\left(m_1 - \overline{m}\right) + D_k\right)\right), \; k = 1,...,4 \qquad (50)$$

Similarly as in the first case, we will substitute equalities (49) and (50) in the associated system (29) and use a method of undefined coefficients. As a result we will receive system of 16 equations with 16 unknowns,

$$\begin{cases}
m_1 A_1 + \overline{m}A_1 + B_1 = m_1 A_1 + m_2 A_2 + m_3 A_3 - m_4 A_4 \\
m_1 B_1 + \overline{m}B_1 = m_1 B_1 + m_2 B_2 + m_3 B_3 - m_4 B_4 \\
m_1 C_1 - \overline{m}C_1 + D_1 = m_1 C_1 + m_2 C_2 + m_3 C_3 - m_4 C_4 \\
m_1 D_1 - \overline{m}D_1 = m_1 D_1 + m_2 D_2 + m_3 D_3 - m_4 D_4 \\
m_1 A_2 + \overline{m}A_2 + B_2 = m_2 A_1 + m_1 A_2 + m_4 A_3 - m_3 A_4 \\
m_1 B_2 + \overline{m}B_2 = m_2 B_1 + m_1 B_2 + m_4 B_3 - m_3 B_4 \\
m_1 C_2 - \overline{m}C_2 + D_2 = m_2 C_1 + m_1 C_2 + m_4 C_3 - m_3 C_4 \\
m_1 D_2 - \overline{m}D_2 = m_2 D_1 + m_1 D_2 + m_4 D_3 - m_3 D_4 \\
m_1 A_3 + \overline{m}A_3 + B_3 = m_3 A_1 - m_4 A_2 + m_1 A_3 + m_2 A_4 \\
m_1 B_3 + \overline{m}B_3 = m_3 B_1 - m_4 B_2 + m_1 B_3 + m_2 B_4 \\
m_1 C_3 - \overline{m}C_3 + D_3 = m_3 C_1 - m_4 C_2 + m_1 C_3 + m_2 C_4 \\
m_1 D_3 - \overline{m}D_3 = m_3 D_1 - m_4 D_2 + m_1 D_3 + m_2 D_4 \\
m_1 A_4 + \overline{m}A_4 + B_4 = m_4 A_1 - m_3 A_2 + m_2 A_3 + m_1 A_4 \\
m_1 B_4 + \overline{m}B_4 = m_4 B_1 - m_3 B_2 + m_2 B_3 + m_1 B_4 \\
m_1 C_4 - \overline{m}C_4 + D_4 = m_4 C_1 - m_3 C_2 + m_2 C_3 + m_1 C_4 \\
m_1 D_4 - \overline{m}D_4 = m_4 D_1 - m_3 D_2 + m_2 D_3 + m_1 D_4
\end{cases},$$

having solve which we will have such common decisions:

$$\begin{cases} x_1 = e^{(m_1+\overline{m})t}\left(-\frac{m_3\overline{m}+m_2m_4}{m_4^{\ 2}-m_3^{\ 2}}A_3+\frac{m_4\overline{m}+m_2m_3}{m_4^{\ 2}-m_3^{\ 2}}A_4\right)+ \\ \quad +e^{(m_1-\overline{m})t}\left(\frac{m_3\overline{m}-m_2m_4}{m_4^{\ 2}-m_3^{\ 2}}C_3+\frac{m_2m_3-m_4\overline{m}}{m_4^{\ 2}-m_3^{\ 2}}C_4\right) \\ x_2 = e^{(m_1+\overline{m})t}\left(-\frac{m_4\overline{m}+m_2m_3}{m_4^{\ 2}-m_3^{\ 2}}A_3+\frac{m_3\overline{m}+m_2m_4}{m_4^{\ 2}-m_3^{\ 2}}A_4\right)+ \\ \quad +e^{(m_1-\overline{m})t}\left(-\frac{m_2m_3-m_4\overline{m}}{m_4^{\ 2}-m_3^{\ 2}}C_3-\frac{m_3\overline{m}-m_2m_4}{m_4^{\ 2}-m_3^{\ 2}}C_4\right) \\ x_3 = e^{(m_1+\overline{m})t}A_3+e^{(m_1-\overline{m})t}C_3 \\ x_4 = e^{(m_1+\overline{m})t}A_4+e^{(m_1-\overline{m})t}C_4 \end{cases} \qquad (51)$$

To find the values of arbitrary constants we use equation (43) and we will obtain

$$\begin{array}{llll} A_1 = \frac{1}{2} & C_1 = \frac{1}{2} & B_1 = 0 & D_1 = 0 \\ A_2 = \frac{1}{2}\frac{m_2}{\overline{m}} & C_2 = -\frac{1}{2}\frac{m_2}{\overline{m}} & B_2 = 0 & D_2 = 0 \\ A_3 = \frac{1}{2}\frac{m_3}{\overline{m}}, & C_3 = -\frac{1}{2}\frac{m_3}{\overline{m}}, & B_3 = 0, & D_3 = 0, \\ A_4 = \frac{1}{2}\frac{m_4}{\overline{m}} & C_4 = -\frac{1}{2}\frac{m_4}{\overline{m}} & B_4 = 0 & D_4 = 0 \end{array}$$

That is, it is possible to write down exponential function of variable from system $D(W, W, 4)$, if to accept $m_i t \approx m_i$

$$Exp(M) = e^{m_1}\left(ch\sqrt{m_2^{\ 2}+m_3^{\ 2}-m_4^{\ 2}}+\frac{sh\sqrt{m_2^{\ 2}+m_3^{\ 2}-m_4^{\ 2}}}{\sqrt{m_2^{\ 2}+m_3^{\ 2}-m_4^{\ 2}}}(m_2e_2+m_3e_3+m_4e_4)\right) \qquad (52)$$

As we consider the case when radical expression is greater than zero, we can rewrite (52) as

$$Exp(M) = e^{m_1}\left(ch\sqrt{\left|m_2^{\ 2}+m_3^{\ 2}-m_4^{\ 2}\right|}+\frac{sh\sqrt{\left|m_2^{\ 2}+m_3^{\ 2}-m_4^{\ 2}\right|}}{\sqrt{\left|m_2^{\ 2}+m_3^{\ 2}-m_4^{\ 2}\right|}}(m_2e_2+m_3e_3+m_4e_4)\right) \qquad (53)$$

3) $m_2^{\ 2}+m_3^{\ 2}-m_4^{\ 2}=0$

In this case, we will have one multiple real root:

$$\lambda_{1,2,3,4} = m_1. \qquad (54)$$

Given (54), the solution of the associated system (29) is looked for in such form:

$$X_k(t) = e^{m_1 t}\left(A_k + B_k t + C_k t^2 + D_k t^3\right), k = 1,...,4 \qquad (55)$$

We calculate the first derivative:

$$\frac{dX_k(t)}{dt} = e^{m_1 t}\left(m_1 A_k + B_k + (m_1 B_k + 2C_k)t + (m_1 C_k + 3D_k)t^2 + m_1 D_k t^3\right),\ k = 1,\ldots,4 \tag{56}$$

Again, similarly as in previous cases, we will substitute equalities (55) and (56) in the associated system (29) and use a method of undefined coefficients. As a result we will receive system of 16 equations with 16 unknowns, having solve which we will have such common decisions.

$$\begin{cases} x_1 = e^{m_1 t}\left(-\frac{m_4}{m_2}A_3 + \frac{m_3}{m_2}A_4 - \frac{m_3}{m_2^{\ 2}}B_3 + \frac{m_4}{m_2^{\ 2}}B_4 + \left(-\frac{m_4}{m_2}B_3 + \frac{m_3}{m_2}B_4\right)t\right) \\ x_2 = e^{m_1 t}\left(-\frac{m_3}{m_2}A_3 + \frac{m_4}{m_2}A_4 - \frac{m_4}{m_2^{\ 2}}B_3 + \frac{m_3}{m_2^{\ 2}}B_4 + \left(-\frac{m_3}{m_2}B_3 + \frac{m_4}{m_2}B_4\right)t\right) \\ x_3 = e^{m_1 t}(A_3 + C_3 t) \\ x_4 = e^{m_1 t}(A_4 + C_4 t) \end{cases} \tag{57}$$

To find the values of arbitrary constants we use equation (43) and we will obtain:

$$\begin{matrix} A_1 = 1 & B_1 = 0 & C_1 = 0 & D_1 = 0 \\ A_2 = 0 & B_2 = m_2 & C_2 = 0 & D_2 = 0 \\ A_3 = 0, & B_3 = m_3, & C_3 = 0, & D_3 = 0. \\ A_4 = 0 & B_4 = m_4 & C_4 = 0 & D_4 = 0 \end{matrix}$$

Thus, the representation of exponent will have the form:

$$Exp(M) = e^{m_1}(e_1 + m_2 e_2 + m_3 e_3 + m_4 e_4) \tag{58}$$

Having compared the received results (47), (53) and (58), we see that representation of exponent for the two first cases are coincides. If to direct on them a radicand of equality (32) to zero, we will receive representation of exponent for the third case.

**Construction of exponent representation of variable from system $D(\boldsymbol{D}, \boldsymbol{D}, 4)$**

In this case we will consider system of numbers as $D(\boldsymbol{D}, \boldsymbol{D}, 4)$. We will build in it exponent representations rom number $M \in D(\boldsymbol{D}, \boldsymbol{D}, 4)$.

Again, the equation (19) can be written down as

$$\frac{dX}{dt} = MX\ . \tag{59}$$

We calculate the right side of equation (59). According to the multiplication table for these numbers we have:

$$M = (m_1 e_1 + m_2 e_2 + m_3 e_3 + m_4 e_4),\ X = (x_1 e_1 + x_2 e_2 + x_3 e_3 + x_4 e_4)$$

$$\begin{aligned} MX &= (m_1 e_1 + m_2 e_2 + m_3 e_3 + m_4 e_4)(x_1 e_1 + x_2 e_2 + x_3 e_3 + x_4 e_4) = \\ &= (m_1 x_1)e_1 + (m_2 x_1 + m_1 x_2)e_2 + \\ &+ (m_3 x_1 + m_1 x_3)e_3 + (m_4 x_1 - m_3 x_2 + m_2 x_3 + m_1 x_4)e_4 \end{aligned}\ . \tag{60}$$

According to (60), equation (59) will have the form:

$$\frac{dX}{dt} = (m_1x_1)e_1 + (m_2x_1 + m_1x_2)e_2 + \\ + (m_3x_1 + m_1x_3)e_3 + (m_4x_1 - m_3x_2 + m_2x_3 + m_1x_4)e_4 . \quad (61)$$

From here it is possible to write down the associated system:

$$\begin{cases} \frac{dx_1}{dt} = m_1x_1 \\ \frac{dx_2}{dt} = m_2x_1 + m_1x_2 \\ \frac{dx_3}{dt} = m_3x_1 + m_1x_3 \\ \frac{dx_4}{dt} = m_4x_1 - m_3x_2 + m_2x_3 + m_1x_4 \end{cases} \quad (62)$$

For solution of the equation (62) it is necessary to find roots of the characteristic equation (24), where

$$M = \begin{pmatrix} m_1 & 0 & 0 & 0 \\ m_2 & m_1 & 0 & 0 \\ m_3 & 0 & m_1 & 0 \\ m_4 & -m_3 & m_2 & m_1 \end{pmatrix},$$ $E$ – single matrix.

Taking into account the aforesaid, the characteristic equation will have the form

$$\begin{vmatrix} m_1 - \lambda & 0 & 0 & 0 \\ m_2 & m_1 - \lambda & 0 & 0 \\ m_3 & 0 & m_1 - \lambda & 0 \\ m_4 & -m_3 & m_2 & m_1 - \lambda \end{vmatrix} = 0 . \quad (63)$$

Having calculated the equation (63), we will receive the following equation:

$$(m_1 - \lambda)^4 = 0 . \quad (64)$$

With (64) it is easy to find characteristic values:

$$\lambda_{1,2,3,4} = m_1 . \quad (65)$$

In this case, common solutions are looked for in following form

$$X_k(t) = e^{m_1 t}\left(A_k + B_k t + C_k t^2 + D_k t^3\right), k = 1,...,4 \quad (66)$$

We calculate the first derivative:

$$\frac{dX_k(t)}{dt} = e^{m_1 t}\left(m_1 A_k + B_k + (m_1 B_k + 2C_k)t + (m_1 C_k + 3D_k)t^2 + m_1 D_k t^3\right), k = 1,...,4 . \quad (67)$$

We will substitute equalities (66) and (67) in the associated system (62) and use a method of undefined coefficients. As a result we will receive system of 16 equations with 16 unknowns, having solve which we will have such common decisions.:

$$\begin{cases} x_1 = e^{m_1 t} \dfrac{1}{m_3} B_3 \\ x_2 = e^{m_1 t} \left( \dfrac{m_2}{m_3} A_3 + \dfrac{m_4}{m_3{}^2} B_3 - \dfrac{1}{m_3} B_4 + \dfrac{m_2}{m_3} B_3 t \right) \\ x_3 = e^{m_1 t} \left( A_3 + C_3 t \right) \\ x_4 = e^{m_1 t} \left( A_4 + C_4 t \right) \end{cases}$$

To find the values of arbitrary constants we use equation (20) and we will obtain:

$$\begin{matrix} A_1 = 1 & B_1 = 0 & C_1 = 0 & D_1 = 0 \\ A_2 = 0 & B_2 = m_2 & C_2 = 0 & D_2 = 0 \\ A_3 = 0, & B_3 = m_3, & C_3 = 0, & D_3 = 0. \\ A_4 = 0 & B_4 = m_4 & C_4 = 0 & D_4 = 0 \end{matrix}$$

That is, it is possible to write down exponential function of variable from system $D(\boldsymbol{D}, \boldsymbol{D}, 4)$, if to accept $m_i t \approx m_i$

$$Exp(M) = e^{m_1} \left( e_1 + m_2 e_2 + m_3 e_3 + m_4 e_4 \right)$$

Exponent representations for all systems, which were considered are given in table 5.

Table 5

**Exponent representations**

| № | HNS | Exponent representation |
|---|---|---|
| 1. | $\boldsymbol{H}$ | $Exp(M)=e^{m_1}\left(\cos\sqrt{m_2{}^2+m_3{}^2+m_4{}^2}\,e_1+\frac{(m_2e_2+m_3e_3+m_4e_4)}{\sqrt{m_2{}^2+m_3{}^2+m_4{}^2}}sh\sqrt{m_2{}^2+m_3{}^2+m_4{}^2}\right)$ |
| 2. | $\boldsymbol{AH}$ | 1) $m_3{}^2+m_4{}^2-m_2{}^2<0$, $m_3{}^2+m_4{}^2-m_2{}^2>0$ <br> $Exp(M)=e^{m_1}\left(ch\sqrt{\lvert -m_2{}^2+m_3{}^2+m_4{}^2\rvert}e_1+\frac{(m_2e_2+m_3e_3+m_4e_4)}{\sqrt{\lvert -m_2{}^2+m_3{}^2+m_4{}^2\rvert}}sh\sqrt{\lvert -m_2{}^2+m_3{}^2+m_4{}^2\rvert}\right)$ <br> 2) $m_3{}^2+m_4{}^2-m_2{}^2=0$ <br> $Exp(M)=e^{m_1}(e_1+m_2e_2+m_3e_3+m_4e_4)$ |
| 3. | $D(\boldsymbol{C},\boldsymbol{D},4)$ | $Exp(M)=e^{m_1}\left(\cos\lvert m_2\rvert e_1+\frac{\sin\lvert m_2\rvert}{\lvert m_2\rvert}(m_2e_2+m_3e_3+m_4e_4)\right)$ |
| 4. | $D(\boldsymbol{W},\boldsymbol{W},4)$ | 1) $m_2{}^2+m_3{}^2-m_4{}^2<0$, $m_2{}^2+m_3{}^2-m_4{}^2>0$ <br> $Exp(M)=e^{m_1}\left(ch\sqrt{\lvert m_2{}^2+m_3{}^2-m_4{}^2\rvert}e_1+\frac{(m_2e_2+m_3e_3+m_4e_4)}{\sqrt{\lvert m_2{}^2+m_3{}^2-m_4{}^2\rvert}}sh\sqrt{\lvert m_2{}^2+m_3{}^2-m_4{}^2\rvert}\right).$ <br> 2) $m_2{}^2+m_3{}^2-m_4{}^2=0$ <br> $Exp(M)=e^{m_1}(e_1+m_2e_2+m_3e_3+m_4e_4)$ |
| 5. | $D(\boldsymbol{D},\boldsymbol{D},4)$ | $Exp(M)=e^{m_1}(e_1+m_2e_2+m_3e_3+m_4e_4)$ |
| 6. | $D(\boldsymbol{W},\boldsymbol{D},4)$ | $Exp(M)=e^{m_1}\left(ch\lvert m_2\rvert e_1+\frac{sh\lvert m_2\rvert}{\lvert m_2\rvert}(m_2e_2+m_3e_3+m_4e_4)\right)$ |

## Conclusions

In this work by means of non-commutative Grassman-Clifford procedure of doubling of systems of second dimension are constructed six various non-commutative HNS of fourth dimension. Their arithmetic and algebraic properties that allows to draw a conclusion of possibility of their using for creation of various mathematical models are investigated.

However, it is necessary to emphasize that we can not assert that these HNS represent different classes of isomorphism, some of them can be isomorphic to each other. This question requires further researches in this area.